\documentclass[11pt]{article}

\usepackage{latexsym,color,graphicx,amsmath,amssymb, xcolor, amsthm,mathtools}

\def\RR{{\mathbb R}}
\def\eps{{\varepsilon}}

\def\p{{\bf p}}
\def\q{{\bf q}}

\newtheorem{thm}{Theorem}
\newtheorem{lem}[thm]{Lemma}

\newtheorem{conj}[thm]{Conjecture}

\theoremstyle{definition}
\newtheorem{defin}[thm]{Definition}
 
\begin{document}

\title{Erd\H{o}s's unit distance problem and rigidity}

\author{
J\'anos Pach\thanks{Alfr\'ed R\'enyi Institute of Mathematics, Budapest, Hungary Email: pach@renyi.hu } \and
Orit E. Raz\thanks{Einstein Institute of Mathematics, The Hebrew University of Jerusalem, Jerusalem, Israel, and Institute for Advanced Study, Princeton NJ, USA Email: oritraz@mail.huji.ac.il} \and
J\'ozsef Solymosi\thanks{University of British Columbia, Vancouver, Canada, and Obuda University, Budapest, Hungary
			Email: solymosi@math.ubc.ca}}
\maketitle

\begin{abstract}
    According to a classical result of Spencer, Szemer\'edi, and Trotter (1984), the maximum number of times the unit distance can occur among $n$ points in the plane is $O(n^{4/3})$. This is far from Erd\H{o}s's lower bound, $n^{1+O(1/\log\log n)}$, which is conjectured to be optimal. We prove a structural result for point sets with nearly $n^{4/3}$ unit distances and use it to reduce the problem to a conjecture on rigid frameworks. This conjecture, if true, would yield the first improvement on the bound of Spencer \emph{et al}. A weaker version of this conjecture has been established by the last two authors.
\end{abstract}

\section{Introduction}

\subsection{Background}
In 1946, Paul Erd\H{o}s~\cite{Erd46, Erd59} raised two different problems for the distribution of distances among $n$ points in the plane or in any fixed metric space $S$:
\smallskip

    \textbf{Problem A.} What is the maximum number of times the same distance can occur among $n$ points in $S$?
\smallskip

    \textbf{Problem B.} What is the minimum number of distinct distances determined by $n$ points of $S$?
\smallskip

Both of these problems remain open today, even on the plane. They generated a lot of research and led to over a thousand papers. This is partially because they turned out to be related to many deep questions in number theory, combinatorics, Fourier analysis, algebraic, incidence, and computational geometry~\cite{GaIS11, Gu16, Ta14, TaV06}. For many related problems, see the monographs~\cite{BrMP05, Ma13, PaA11}.

The two problems are obviously interconnected. If we obtain an upper bound $u(n)$ for the quantity in Problem A, it gives the lower bound $\binom{n}{2}/u(n)$ for the number of distinct distances in Problem B. In particular, Erd\H os conjectured that in the plane $u(n)=n^{1+O(1/\log\log n)},$ which is attained by a $\sqrt{n}\times\sqrt{n}$ piece of the square grid. 
A matching upper bound is known only in exceptional cases where the vectors determining the unit distances are from a restricted family \cite{SoScZe,Sch}.

Unfortunately, despite many efforts, the best known upper bound for the number of times that (say) the \emph{unit distance} can occur among $n$ points in the plane is 
\begin{equation}
u(n)=O(n^{4/3}).
\end{equation}
This was first proved by Spencer, Szemer\'edi, and Trotter~\cite{SpSzT84}. Since then, the problem has been tackled by Clarkson, Edelsbrunner, Guibas, Sharir, and Welzl~\cite{ClEG90}, Sz\'ekely~\cite{Sz97}, Pach and Tardos~\cite{PaT06}, using different approaches, all yielding precisely the same upper bound, $O(n^{4/3})$. None of these approaches offered much hope for possible improvement. 

To partially explain why all previous methods yielded the same upper bound, we reformulate the unit distance problem as an incidence problem between points and unit circles in $\RR^2$:
Let $P$ be a set of $n$ points in $\RR^2$, and let $C$ be a set of $n$ unit circles in $\RR^2$. Consider the set of point-circle incidences
$$
I(P,C)=\{(p,c)\in P\times C\mid p\in c\}.
$$
Observe that 
$$
\max_{|P|=n,|C|=n} I(P,C)=\Theta(u(n)).
$$
That is, the problem of upper bounding $u(n)$ is equivalent to the problem of upper bounding the number of incidences between $n$ points and $n$ unit circles in $\RR^2$.

Bounding the number of incidences between points and curves in $\RR^d$ is a well-studied problem in combinatorial geometry. Specifically, the instance of point-line incidences in $\RR^2$ is the classical result of Szemer\'edi and Trotter from 1984:
\begin{thm}[{\bf Szemer\'edi--Trotter~\cite{SzeTro}}]\label{SzT}
    Let $P$ be a set of $m$ points in $\RR^2$, and let $L$ be a set of $n$ lines in $\RR^2$. Then 
$$
I(P,L)=O(m^{2/3}n^{2/3}+m+n).
$$
Moreover, this bound is tight.
\end{thm}
In the special case  $m=n$, Theorem~\ref{SzT} yields $I(P,L)=O(n^{4/3})$, and this bound is best possible. If we want to obtain a better bound for the incidence problem between points and unit circles in the plane, one has to come up with a method that does not apply to point-line incidences; otherwise, we are doomed to fail. This is ``the" reason why the unit distance problem appears to be so difficult. All known methods, with the exception of the one developed in \cite{PaT06}, also apply to the point-line incidence problem, hence, they cannot yield any bound better than $O(n^{4/3})$. The approach in \cite{PaT06} can be used to bound the number of incidences between points and some other families of curves, for which the bound is tight; see \cite{Val,SoSza}.

As for Erd\H{o}s's other problem (Problem B), concerning the smallest number of distinct distances, the lower bound has been steadily improved over the years by Moser~\cite{Mo52}, Chung, Szemer\'edi, and Trotter~\cite{ChSzT92}, Solymosi and Cs.~T\'oth~\cite{SoT01}, and Katz and Tardos~\cite{KaT04}. 
In May 2010, Elekes and Sharir reduced the question to an incidence problem in $\RR^3$. A couple of months later, in November of the same year, Guth and Katz~\cite{GuK15} achieved a major breakthrough by solving the incidence problem of Elekes and Sharir. They deduced that any set of $n$ points in the plane determines at least $cn/\log n$ distinct distances, where $c>0$ is a suitable constant. This is only a factor of $\sqrt{\log n}$ smaller than the conjectured minimum, which is attained again by the square grid construction.

\subsection{Graph rigidity}
Before formulating our results, we need to recall some basic notions from graph rigidity theory.  

We review some standard definitions from rigidity theory. For more details, see, e.g., Asimow and Roth~\cite{AsimowRoth1}. Let \(G = (V,E)\) be a graph. A \emph{realization} \(\p\) of $G$ in \( \mathbb{R}^2\) is an embedding (not necessarily injective) of the vertex set \(V = \{1,\ldots,|V|\}\) of $G$ in $\RR^2$. That is, \[\p = (p_1,..,p_{|V|}) \in \mathbb{R}^2 \times ... \times \mathbb{R}^2 \cong \mathbb{R}^{2|V|}.\] 
A pair $(G,\p)$ of a graph $G$ and a realization $\p$ is called a \emph{framework}. 

Define the \emph{edge function} of $G$
to be the map \(f_G: \mathbb{R}^{2|V|} \to \mathbb{R}^{|E|}\) given by \[(p_1,\ldots,p_{|V|})\mapsto \left(||p_i - p_j||^2\right)_{\{i,j\}\in E},\]
where we fix an arbitrary order on the edges of $G$.
Note that \(f_G\) is a polynomial map.

For a given graph $G$, let $\p$ and $\q$ be a pair of realizations of $G$ in $\RR^2$. We say that the corresponding frameworks, $(G,\p)$ and $(G,\q)$, are \emph{equivalent} if \(f_{G}(\p) = f_{G}(\q)\), that is, if $\|p_i-p_j\|=\|q_i-q_j\|$ holds for every edge $\{i,j\}\in E$. We say that the frameworks $(G,\p)$ and $(G,\q)$ are \emph{congruent} if $\|p_i-p_j\|=\|q_i-q_j\|$ for every size-2 subset $\{i,j\}\subset V$ (here $\{i,j\}$ is not necessarily an edge in $G$). Equivalently, $(G,\p)$ and $(G,\q)$ are congruent if there exists an isometry $R$ of $\RR^2$ such that $R(p_i)=q_i$, for every $i\in V$.

\begin{defin}[{\bf  Framework Rigidity}]\label{def:frameworkrigidity}
We say that \((G,\p)\) is a \emph{rigid framework} if there exists a neighborhood $U\subset (\RR^2)^{|V|}$ of $\p$, such that for every $\q\in U$, if $(G,\p)$ and $(G,\q)$ are equivalent, then they are necessarily also congruent. Equivalently, if there exists a neighborhood $U$ of $\p$ such that
$$
f_G^{-1}(f_G(\p))\cap U=f_{K_{|V|}}^{-1}(f_{K_{|V|}}(\p))\cap U,$$
where $K_{|V|}$ is the complete graph on $|V|$ vertices.
\end{defin}

Note that the notion of framework rigidity depends on both the graph $G$ and the realization $\p$. That is, for a given graph $G$, and realizations $\p$ and $\q$, it might happen that $(G,\p)$ is rigid and $(G,\q)$ is not rigid. 

A realization of $G$ is called \emph{generic} if the set of all coordinates of all of its points is algebraically independent over the rationals. It turns out that all {\it generic} realizations of $G$ behave the same. 
That is, if $G$ is fixed and $\p$ and $\q$ are some generic realizations of $G$, then $(G,\p)$ is rigid if and only if $(G,\q)$ is rigid. In this sense, rigidity is a property of the graph $G$, independent of the specific generic realization $\p$ one considers.  
      
\begin{defin}[{\bf Graph Rigidity}]\label{def:graphrigidity}
We say that $G=([n],E)$ is a \emph{rigid graph} in $\RR^2$, if $(G,\p)$ is a rigid framework, for every generic realization  $\p\in(\RR^2)^n.$
\end{defin}

Most results in rigidity theory are concerned with the notion of generic rigidity, as in Definition~\ref{def:graphrigidity}. In the present note, we study realizations $\p$ arising from configurations of $n$ points in the plane that maximize the number of unit distances. Such configurations are highly \emph{non-generic}. Hence, for us \emph{framework rigidity}, as in Definition~\ref{def:frameworkrigidity}, will be more relevant, and we will focus on this notion.   

Specifically, given a framework $(G,\p)$, we aim to find a rigid sub-framework within it. Note that this is an easy task if we restrict our attention to \emph{generic} realizations. It is not hard to show that every sufficiently dense graph contains a \emph{generically} rigid subgraph (see, e.g., \cite{Con} for technical details).

\begin{lem} 
    Let $G=(V,E)$, let $\eps>0$ be fixed, and suppose that $|E|\ge n^{1+\eps}$.    
    Then there exists $G'\subset G$ with $|V(G')|\ge 3$ such that $G'$ is generically rigid in $\RR^2$.
\end{lem}

However, if our graph is given together with a realization, there is not much theory about this problem. The only useful result in this direction was obtained by the last two authors in \cite{RazSol}.
\begin{thm}[{\bf Raz and Solymosi \cite{RazSol}}]\label{rigidsubframework} Let $G=([n],E)$, let $\alpha>1/2$, and suppose that 
\begin{equation}\label{edges}
|E|=\Omega(n^{1+\alpha}).
\end{equation} 
Let $\p\in(\RR^2)^{n}$ be a realization of $G$ with the property that for every vertex $v\in [n]$, the neighbours of $v$ are not embedded into a common line.  
Then there exists a subgraph $G'\subset G$, with $|V(G')|\ge 4$, such that $(G', \p|_{V(G')})$ is a rigid framework, provided $n$ is large enough.
\end{thm}

The smallest value of $\alpha$ for which the conclusion of Theorem~\ref{rigidsubframework} holds is not known. As was observed in \cite{RazSol}, we need to assume at least $|E|=\Omega( n\log n)$; otherwise, the statement is false.

    In \cite{RazSol}, Theorem~\ref{rigidsubframework} was proved under the slightly stronger assumption that there are no \emph{three} vertices of $G$ that are mapped by $\p\in(\RR^2)^{n}$ into collinear points. However, essentially the same proof also gives Theorem~\ref{rigidsubframework}.

\subsection{Main results}

Our goal is to suggest a new approach to tackle the unit distance problem (Problem A) by reducing it to a rigidity problem.

For any point set $P$ in the plane, let $u(P)$ denote the number of unit distance pairs determined by $P$, that is, the number of unordered pairs $\{p,q\}\subseteq P$ with $\|p-q\|=1$. Given a graph $G=(V,E)$ and a realzation $\p$ of its vertex set into $\RR^{2|V|}$ such that, $\p$ is injective, and the endpoints of every edge are mapped into two points whose distance is $1$, we call $\p$ a \emph{unit embedding} of $G$.

In the sequel, when we compare two functions, we will often use the notation $f(n)\lesssim g(n)$ instead of $f(n)=O(g(n))$. If we have $f(n)\lesssim g(n)$ and $g(n)\lesssim f(n),$ then we write $g(n)\approx h(n)$.


We tacitly assume that there exist $n$-element point sets $P$ with $u(P)$ close to the currently best known upper bound $n^{4/3}$, and we study their structure.
Our main technical result is the following.

\begin{thm}[{\bf Structure Theorem}]\label{main:structure}
Let $h(n)$ be a function tending to $\infty$, as $n\rightarrow\infty$. Let $P$ be a set of $n$ points in $ \RR^2$ satisfying $u(P)\ge\frac{n^{4/3}}{h(n)}$, and suppose that $n$ is sufficiently large.

 Then there exist
\begin{enumerate}
    \item a subset $P'\subset P$ with $|P'|\approx n^{1/3}h(n)^4$;
    \item bipartite graphs $G_i=(U_i\cup V_i, E_i)$ for every $i\; (1\le i\le k),$ where $k\gtrsim n^{2/3}/h(n)^{5}$,  such that $2\le |U_i|,|V_i|\lesssim h(n)^6$
    and  $|E_i|\gtrsim h(n)^{7};$
    \item unit embeddings $\p^{(i)}$ of $G_i$ into $(\RR^2)^{|U_i|+|V_i|}$ for every $i\; (1\le i\le k),$ such that
    
    \indent (i) the sets $\p^{(1)}(U_1),\ldots, \p^{(k)}(U_k)$ are pairwise disjoint subsets of $P$,\\
    \indent (ii) 
    the sets $\p^{(1)}(V_1),\ldots, \p^{(k)}(V_k)$ are subsets of $P'$.
\end{enumerate}
    \end{thm}
    
We conjecture that the statement of Theorem~\ref{rigidsubframework} remains true for any $\alpha\ge 1/6$. In other words, we state the following conjecture.

\begin{conj}[{\bf Rigidity Conjecture}]\label{rigidityconj}
Let $G=([n],E)$ be a graph with $|E|\gtrsim n^{7/6}.$ 
Let $\p\in(\RR^2)^{n}$ be a realization of $G$ with the property that for every vertex $v\in [n]$, the neighbours of $v$ are not embedded into a common line.  
Then there exists a subgraph $G'\subset G$ with $|V(G')|\ge 4$ such that $(G', \p|_{V(G')})$ is a rigid framework.
\end{conj}

Provided that Conjecture~\ref{rigidityconj} is true, our approach would yield the first improvement of the classical upper bound, $O(n^{4/3})$, of Spencer, Szemer\'edi, and Trotter~\cite{SpSzT84} on the number of unit distances for more than 40 years. Notably, we will deduce the following result.

\begin{thm}\label{mainthm}
Conjecture~\ref{rigidityconj} implies that the maximum number of unit distance pairs, $u(n)$, determined by $n$ points in the plane satisfies
$$u(n)=O\left(\frac{n^{4/3}}{\log^{1/12} n}\right).$$
\end{thm}

The proofs of Theorems~\ref{mainthm} and~\ref{main:structure} are given in Sections~\ref{sec:proofofmain} and~\ref{sec:proofofstructure}, respectively.




\section{Proof of Theorem~\ref{mainthm}}\label{sec:proofofmain}
    Let $P\subset \RR^2$ with $|P|=n$, and assume for contradiction that $u(P)\ge \frac{n^{4/3}}{h(n)}$, for some function $h=h(n)$ that tends to zero as $n$ tends to infinity.
    Let $P'\subset P$, $k$, $G_1=(U_1\cup V_1,E_1),\ldots, G_k=(U_k\cup V_k, E_k)$ and $\p^{(1)},\ldots,\p^{(k)}$ be as given by Theorem~\ref{main:structure}. 

    Note that, since the realizations $\p^{(i)}$ are unit embeddings, the assumption in Conjecture~\ref{rigidityconj}, that the neighbours of a vertex $v$ are not embedded to a common line,  applies automatically. Thus, assuming that Conjecture~\ref{rigidityconj} is true, it follows that, for every $i=1,\ldots, k$, there exists a subgraph $G_i'=(U_i'\cup V_i',E_i')\subset G_i$ such that $G_i'$ has at least 4 vertices and  the unit framework $(G_i',\p^{(i)}|_{U_i'\cup V_i'})$ is rigid. Since $G_i'$ is bipartite, this implies, in particular, that $|U_i|,|V_i|\ge 2$.
    
    Recall that, by Theorem~\ref{main:structure} item 2, we have that for each $i$,  $2\le |U_i|,|V_i|\le h^6$. Thus also $2\le |U_i'|,|V_i'|\le h^6$. By the pigeonhole principle, there exists some bipartite graph $H=(U\cup V, E)$, with $2\le |U|,|V|\le h^6$, such that for at least 
    $$
    k/2^{h^{12}},$$
indices $i$, we have that $G_i'$ is isomorphic to $H$. Note that a rigid framework in the plane, on $m$ vertices,  has at most  $9^m$ distinct non-congruent embeddings that induce the same edge lengths (see Milnor  \cite[Theorem 2]{Milnor}). Applying the pigeonhole principle once again, we conclude that there are at least 
 $$
    k/(2^{h^{12}}9^{2h^6})$$ 
    indices $i$, for which the frameworks $(G_i',\p^{(i)}|_{U_i'\cup V_i'})$ are pairwise congruent. Let $I$ denote the subset of such indices. So $|I|\ge k/(2^{h^{12}}9^{2h^6})$.

    Let $v_1,v_2\in V$. Then there exists a number $a>0$ such that, for every $i\in I$, the embedding $\p^{(i)}$ maps $v_1,v_2$ to a pair of points in $P'$ that are at distance $a$ from each other.

    We claim that for every $p,q\in P'$, there are at most two indices $i\in I$ such that  $v_1$ is embedded by $\p^{(i)}$ to $p$ and $v_2$ is embedded by $\p^{(i)}$ to $q$. 
Indeed, note that by fixing the embedding of $v_1,v_2$, we have that $\p^{(i)}$ must be one of at most two possible realizations of $H$ (which can be obtained from each other by a reflection through the line $pq$). Thus, the embedding of $U_i'$ is determined up to two possibilities. Since the sets $\p^{(i)}(U_i)$ are pairwise disjoint, it follows that $i$ must be one of at most two possible indices in $I$. 
    Similarly,  there
    are at most two indices $i\in I$ such that  $\p^{(i)}$ embeds $v_1$ to $q$ and $v_2$ to $p$. 

    We conclude that the number of  pairs in $P'$ determining distance $a$ is at least
 $$
    k/(4\cdot 2^{h^{12}}\cdot 9^{2h^6}).$$ 
Since the same distance can be repeated at most $\lesssim |P'|^{4/3}$ times in $|P'|$, we have
$$
    k/(4\cdot 2^{h^{12}}\cdot 9^{2h^6})\lesssim |P'|^{4/3}\approx (n^{1/3}h^4)^{4/3}.
$$
Using the lower bound on $k$, we obtain that
\begin{align*}
    \frac{n^{2/3}}{h^52^{h^{12}+2h^6+1}}\lesssim n^{4/9}h^{16/3}
\end{align*}
or
\begin{align*}
    n^{2/9}\lesssim 2^{h^{12}+2\log(9)h^6+1+(31/3)\log h}. 
\end{align*}
This yields a contradiction if $h(n)\approx \log^{1/12}n$, completing the proof of the theorem.\hfill\qed

\section{Proof of
Theorem~\ref{main:structure}\label{sec:proofofstructure}
}
For the proof, we need the following result of Guth~\cite{Guth}.
\begin{thm}[{{\bf Guth~\cite[Theorem 0.3]{Guth}}}]\label{Guthpartition}
Let $\Gamma$ be a set of $k$-dimensional varieties in $\RR^d$, 
each defined by at most $b$ polynomial equations of degree at most $\delta$.
Then, for any $D\ge 1$, there is a nonzero polynomial $f$ of degree at most $D$ such that each connected component of $\RR^d\setminus \{f=0\}$
intersects $\lesssim D^{k-d}|\Gamma|$ varieties $\gamma\in \Gamma$, where the constant of proportionality depends on $d,\delta$, and $b$. 
\end{thm}

Let $P\subset \RR^2$ with $|P|=n$, where $n$ is sufficiently large. Assume that $$u(P)\ge \frac{n^{4/3}}{h(n)},$$    
for some  function $h=h(n),$ which tends to $\infty$ as $n\rightarrow\infty$.

Let $C$ denote the set of unit circles centred at the points of $P$. By our assumption on $P$, we have that the number of incidences between $P$ and $C$ is 
$$
|I(P,C)|\ge \frac{n^{4/3}}{h} 
$$
\paragraph{First partitioning.}
Set 
$$r=\frac{n^{1/3}}{h^2}.$$ 
By Theorem~\ref{Guthpartition}, there exists $f\in\RR[x,y]$, with $\deg(f)\le r$, whose zero set partitions the plane into $\lesssim r^2$ cells such that each cell contains at most $\lesssim n/r^2$ points of $P$ and meets at most $\lesssim n/r$ unit circles of $C$.

The number of incidences that occur on the zero set of $f$ is at most
$$
|I_0(P,C)|\le  nr+n+r^2<\frac12 |I(P,C)|,
$$
if $n\ge n_0$ is sufficiently large.

Thus, in what follows, we assume without loss of generality that $P\subset \RR^2\setminus Z(f)$ and that no circle in $C$ is contained in $Z(f)$.

Let $\Omega$ denote the set of open connected components of $\RR^2\setminus Z(f)$. For $\omega\in \Omega$, let $P_\omega$ denote the set of points in $P$ contained in $\omega$, and let $C_\omega$ denote the set of circles in $C$ that intersect $\omega$. 
By our choice of $f$, we have 
$$
|\Omega|\lesssim \frac{n^{2/3}}{h^4},$$
and, for every $\omega\in \Omega$, 
$$
|P_\omega|\lesssim n/r^2=n^{1/3}h^4
$$
and 
$$
|C_\omega|\lesssim n/r=n^{2/3} h^2.
$$

Recall our assumption that all incidences in $I(P,C)$ occur within the cells. Then, by the pigeonhole principle, there exists a cell $\omega_0\in \Omega$ such that 
$$
|I(P_{\omega_0},C_{\omega_0})|\gtrsim \frac{n^{4/3}}{h}/\frac{n^{2/3}}{h^4}
=n^{2/3}h^3.
$$

Let
$Q$ denote the set of centers of the unit circles in $C_{\omega_0}$. Let $D$ denote the set of unit circles centred at the points of $P_{\omega_0}$. 
We have 
\begin{align*}
|Q|&\lesssim n^{2/3}h^2,\\
|D|&\lesssim n^{1/3}h^4,\quad\text{and}\\
|I(Q,D)|&\gtrsim n^{2/3}h^3.
\end{align*}

\paragraph{Second partitioning.}
Set 
$$t \approx \frac{n^{1/3}}{h^2}.$$
Let $g\in\RR[x,y]$ be a bivariate polynomial of degree at most $t$ such that $\RR^2\setminus Z(g)$ partitions $\RR^2$ into at most $\lesssim t^2$ connected components, each  of which contains at most 
$\lesssim |Q|/t^2\lesssim h^6$ 
points of $Q$, and 
meets 
at most 
$\lesssim |D|/t\lesssim h^6$ 
unit circles in $D$. 

Note that the number of incidences in $I(Q,D)$ that occur on the zero set of $g$ is bounded by 
$$
|I_0(Q,D)|\le  t|D|+|Q|+t^2<\frac12|I(Q,D)|,
$$
for $n\ge n_0$ sufficiently large. 
Thus, without loss of generality, we may assume that all incidences in $I(Q,D)$ occur within the cells.

Let $\Pi$ denote the open connected components of $\RR^2\setminus Z(g)$. For $\pi\in \Pi$, let $Q_\pi$ denote the set of points of $Q$ contained in $\pi$, and let $D_\pi$ denote the set of circles in $D$ that meet $\pi$. 
By our choice of $g$, we get
$$
|\Pi|\lesssim \frac{n^{2/3}}{h^4}$$
and, for every $\pi\in \Pi$, we have
\begin{align*}
|Q_\pi|&\lesssim h^6\quad{and}\\
|D_\pi|&\lesssim h^6,
\end{align*}
as was already noted above.

Note that cells $\pi$ with at most $\le ch^7$ incidences, for some constant $c>0$, contribute a total of at most 
$$
\lesssim \frac{n^{2/3}}{h^4}\cdot c h^7= cn^{2/3}h^3<\frac12 I(Q,D)$$
incidences, provided that $c>0$ is chosen sufficiently small.

Let $\Pi'\subset \Pi$ be the subset of cells $\pi$ for which $|I(Q_\pi,D_\pi)|\ge ch^7$. Observe that 
\begin{equation}\label{lowerk}|\Pi'| \gtrsim \frac{n^{2/3}}{h^{5}}.
\end{equation}
Indeed, by the Szemer\'edi--Trotter bound, in each cell $\pi\in \Pi$, we have
$|I(Q_\pi,D_\pi)|\lesssim (|Q_\pi||D_\pi|)^{2/3}\approx h^8$.
Thus, we have 
$$
|\Pi'|\cdot h^8
\gtrsim \frac12 |I(Q,D)|\gtrsim 
n^{2/3}h^3,$$
which 
implies \eqref{lowerk}.

Finally, let $P':=P_{\omega_0}$, or, equivalently, the centers of the circles in $D$. Note that $P'\subset P$.
For every $\pi\in \Pi'$, let $G_\pi$ denote the incidence graph between the points in $Q_\pi$ and the unit circles in $D_\pi$. 
Observe that $G_\pi$ is a bipartite graph with parts $Q_\pi, D_\pi$, each of cardinality at most $\lesssim h^6$, and that the number of edges in  $G_\pi$ is $|E(G_\pi)|\gtrsim h^7$. Moreover, identifying the circles in $D_{\pi}$ with their centers, we get that $(Q_\pi, D_\pi)$ is a \emph{unit embedding} of $G_\pi$, by construction.
Observing also that the sets $Q_\pi\subset P$ are pairwise disjoint (by the properties of the partition induced by the partitioning polynomial $g$), and in view of \eqref{lowerk}, the proof of Theorem~\ref{main:structure} is complete.\hfill\qed


\section{Acknowledgements}
The authors thank Zvi Shem-Tov, Joshua Zahl, and Frank de Zeeuw for their valuable discussions and assistance in shaping the paper. 

The work of the first-named author was supported by ERC Advanced Grant no. 882971  ``GeoScape'' and by the National Research, Development and Innovation Office NKFIH Grant no. K-131529.
The research of the second-named author was supported by the Charles Simonyi Endowment.  
The third-named author's research was supported by an NSERC Discovery grant and by the National Research Development and Innovation Office of
Hungary, NKFIH, Grants no. KKP133819 and Excellence 151341.


\begin{thebibliography}{9}

\bibitem{AsimowRoth1}
Leonard Asimow and Ben Roth. 
The rigidity of graphs. 
{\it Trans. Amer. Math. Soc.} 
245 (1978), 279--289.



\bibitem{AsimowRoth2}
Leonard Asimow and Ben Roth. 
The rigidity of graphs II. 
{\it J. Math. Anal. Appl.}, 68 (1979), 171--190.

\bibitem{BorStr}
C. Borcea and I. Streinu,
The number of embeddings of minimally rigid graphs,
\emph{Discrete Comput. Geom.} 31 (2004), 287--303, DOI: 10.1007/s00454-003-2902-0.


\bibitem{BrMP05}
Peter Brass, William O.J. Moser, and J\'anos Pach. \emph{Research Problems in Discrete Geometry.} New York, Springer, 2005.



\bibitem{Con}
Robert Connelly. Generic global rigidity. \emph{Discrete \& Computational Geometry} \textbf{33} (2005), 549--563. 

\bibitem{ChSzT92}
Fan Chung, Endre Szemerédi, and William T. Trotter. The number of different distances determined by a set of points in the Euclidean plane. \emph{Discrete \& Computational Geometry} \textbf{7}, no. 1 (1992), 1--11.

\bibitem{ClEG90} Kenneth L. Clarkson, Herbert Edelsbrunner, Leonidas J. Guibas, Micha Sharir, and Emo Welzl. Combinatorial complexity bounds for arrangements of curves and spheres. \emph{Discrete \& Computational Geometry} \textbf{5}, no. 2 (1990), 99--160.

\bibitem{Erd46}
Paul Erd\H{o}s. On sets of distances of $n$ points.
\emph{American Mathematical Monthly} \textbf{53}, no. 5 (1946), 248--250.

\bibitem{Erd59}
Paul Erd\H{o}s. On sets of distances of $n$ points in Euclidean space. \emph{Magyar Tudom\'anyos Akad\'emia Matematikai Kutat\'o Int\'ezet\'enek K\"ozlem\'enyei} \textbf{5}, no. 1--2 (1959), 165--169.

\bibitem{GaIS11}
Julia Garibaldi, Alex Iosevich, and Steven Senger. \emph{The Erdos Distance Problem.} American Mathematical Soc., 2011.

\bibitem{Guth}
Larry Guth. Polynomial partitioning for a set of varieties.
\emph{Mathematical Proceedings Cambridge Philosophical Society} \textbf{159} (2015), 459--469.

\bibitem{Gu16}
Larry Guth. \emph{Polynomial Methods in Combinatorics.} American Mathematical Society, 2016.

\bibitem{GuK15}
Larry Guth and Nets Hawk Katz. On the Erd\H{o}s distinct distances problem in the plane. \emph{Annals of Mathematics} (2015), 155--190.

\bibitem{KaT04}
Nets Hawk Katz and G\'abor Tardos. A new entropy inequality for the Erd\H{o}s distance problem. \emph{Contemporary Mathematics}, Vol. \textbf{342} (2004), 119--126.

\bibitem{Ma13}
Ji\v ri Matou\v sek. \emph{Lectures on Discrete Geometry.} Springer Science \& Business Media, 2013.

\bibitem{Milnor}
J. Milnor, On the Betti numbers of real varieties, \emph{Proc. AMS} 15 (1964), 275--280

\bibitem{Mo52}
Leo Moser. On the different distances determined by $n$ points. \emph{American Mathematical Monthly} \textbf{59}, no. 2 (1952), 85--91.

\bibitem{PaA11}
J\'anos Pach and Pankaj K. Agarwal. \emph{Combinatorial Geometry.} John Wiley \& Sons, 2011.

\bibitem{PaT06}
J\'anos Pach and G\'abor Tardos. Forbidden paths and cycles in ordered graphs and matrices. \emph{Israel Journal of Mathematics} \textbf{155}, no. 1 (2006), 359--380.

\bibitem{RazSol} Orit E. Raz and J\'ozsef Solymosi. Dense Graphs Have Rigid Parts. \emph{Discrete \& Computational Geometry} \textbf{69} (2023), 1079–1094. 

\bibitem{SoScZe}
{Ryan Schwartz, J\'ozsef Solymosi, and Frank de Zeeuw}. {Rational distances with rational angles}. \emph{Mathematika}, \textbf{58}, {2}, (2012), {409--418}.

\bibitem{Sch}
Ryan Schwartz. {Using the subspace theorem to bound unit distances}, \emph{Mosc. J. Comb. Number Theory}, \textbf{3}, {1}, (2013), {108--117}.

\bibitem{SoSza}
J\'ozsef Solymosi and Endre Szabó.
Classification of maps sending lines into translates of a curve,
\emph{Linear Algebra and its Applications},
Volume 668, 2023, 161-172,



\bibitem{SoT01}
J\'ozsef Solymosi and Csaba D. T\'oth. Distinct distances in the plane. \emph{Discrete \& Computational Geometry} \textbf{25} (2001), 629--634.


\bibitem{SpSzT84}
Joel Spencer, Endre Szemer\'edi, and William T. Trotter. Unit distances in the Euclidean plane. In: \emph{Graph Theory and Combinatorics}, Academic Press, 1984, 294--304.

\bibitem{Sz97}
L\'aszl\'o Sz\'ekely. Crossing numbers and hard Erd\H os problems in discrete geometry. \emph{Combinatorics, Probability and Computing} \textbf{6}, no. 3 (1997), 353--358.

\bibitem{SzeTro}
Endre Szemer\'edi and William T. Trotter.
Extremal problems in discrete geomtery,
\emph{Combinatorica}, \textbf{3} (1983), 381--392.

\bibitem{Ta14}
Terence Tao. Algebraic combinatorial geometry: the polynomial method in arithmetic combinatorics, incidence combinatorics, and number theory. In: \emph{EMS Surveys in Mathematical Sciences} \textbf{1}, no. 1 (2014), 1--46.

\bibitem{TaV06}
Terence Tao and Van H. Vu. \emph{Additive Combinatorics.} Cambridge University Press, 2006.

\bibitem{Val}
P. Valtr. Strictly convex norms allowing many unit distances and related touching questions, manuscript 2005.

\end{thebibliography}
\end{document}